\documentclass[12pt,a4paper]{article}
\usepackage[T1]{fontenc}

\usepackage{amsmath,amssymb,amsfonts}        
\usepackage{mathptmx}
\usepackage[scaled=0.91]{helvet}
\usepackage{courier}
\usepackage[dvips]{graphicx}
\usepackage{psfrag}
\usepackage{color}
\def\ds{\displaystyle}
\def\e{\epsilon}
\def\RR{\hbox{I\kern-.2em\hbox{R}}}
\def\NN{\hbox{I\kern-.2em\hbox{N}}}
\newcommand{\norme}[2]{\ensuremath{\| #1 \|_{#2}}}
\providecommand{\myceil}[1]{\ensuremath{\left \lceil #1 \right \rceil }}
\title{A Front Fixing Implicit Finite Difference Method for the American Put Options Model}
\author{Riccardo Fazio\footnote{Corresponding author's e-mail: rfazio@unime.it}, Alessandra Insana and Alessandra Jannelli\\
Department of Mathematics and Computer Science, \\
Physical Sciences and Earth Sciences\\
Viale F. Stagno d'Alcontres 31, I--98166 Messina, Italy}

\date{\today}

\begin{document}
\maketitle

\begin{abstract}
In this paper, we present an implicit finite difference method for the numerical solution of the Black-Scholes model
of American put options without dividend payments. 
We combine the proposed numerical method by using a front fixing approach
where the option price and the early exercise boundary are computed simultaneously.
Consistency and stability properties of the method are studied. 
We choose to improve the accuracy of the computed solution via a mesh refinement based on Richardson's extrapolation.
Comparisons with some proposed methods for the American options problem are carried out to validate 
the obtained numerical results and to show the efficiency of the proposed numerical methods.
Finally, by \textit{a posteriori} error estimator, we find a suitable computational grid requiring that the computed solution
verifies a prefixed tolerance.
\end{abstract}

\smallskip

\smallskip


\section{Introduction}
American options are contracts allowing the holder to exercise the option, to sell or to buy an asset, at a certain price at any
time prior to and including its maturity date. 
The pricing of American options plays an important role both in theory and in the real derivative markets.
The option pricing model, developed by Black and Scholes  \cite{Black:POC:1973} and extended by Merton \cite{Merton:1973:TRO}, gives rise to a partial differential equation governing the value of an option.
Schwartz \cite{Schwartz:1977;VWI}, Brennan and Schwartz \cite{Brennan:1977:VAP,Brennan:1978:FDM} were the first to apply a finite difference method to price American options. 
The accuracy of their finite difference method was proved by Jaillet et al. \cite{Jaillet:1990:VIP}. 
Other papers present finite difference methods, (see, for example, Hull and White \cite{Hull:1990:VDS}, Duffy \cite{Duffy:2006:FDM}, Wilmott et al. \cite{Wilmott:1995:MFD}). 
An alternative approach is based on front-tracking methods that keep track of the boundary and discretize the problem in changing domain;
they have been considered in \cite{Han:2003:FNM}, \cite{Pantazopoulos:1998:FFD}, \cite{Tangman:1977:FHF}. 

When we price an American option we also need to determine the optimal exercise moment as a function of the value of the underlying asset. 
This leads to formulating free boundary of the non-linear problem for the price of the American option looking for a boundary changing in time to maturity, known as the optimal exercise boundary.
In particular, the American call option problem is a free boundary problem defined on a finite domain.
On the other hand, the American put option problem is a free boundary problem defined on a semi-infinite domain so that it is a non-linear model complicated by a boundary 
condition at infinity.
The difficulty associated with a free boundary problem can be reduced by using a front-fixing method, an approach relying on a change of variable to map the changing domain into a fixed domain. 

The front-fixing approach has been considered in several papers.
In Holmes and Yang \cite{Holmes:2008:FFE}, a front-fixing finite element method was used. 
Tangman et al. \cite{Tangman:1977:FHF} introduced a fourth-order accurate finite difference scheme. 
In \cite{Wu:1997:FFM}, the original problem is transformed into one more manageable equation with coefficients not depending on the spatial variable and an explicit update for the location of the free boundary at each time step was used, and in Zhu et al. \cite{Zhu:2003:ASS}, the secant method was employed to solve the non-linear problems. 
Moreover, in Tangman et al. \cite{Tangman:1977:FHF} and \cite{Zhu:2003:ASS}, the difference between the prices of American options and European options was computed.
{S}ev\v{c}ovi\v{c} \cite{Sevcovic:2007:OEB} studied and approximated the early exercise boundary for a
class of nonlinear Black-Scholes equations, applying a fixed-domain transformation and using an operator splitting iterative numerical technique.
Lauko and \v{S}ev\v{c}ovi\v{c} \cite{lauko2010comparison} introduced a local iterative numerical scheme and compared analytical and numerical solution for the early exercise boundary position computation of the American put option.
In \cite{Zhang:2009:HFD}, Zhang and Zhu presented a predictor-corrector method after the fixed-domain transformation.
Nielsen et al. \cite{Nielsen:2002:PFF} used both explicit and implicit schemes for solving Black-Scholes model of American
options. 
While Company et al. \cite{Company:2014:SAP} presented an explicit finite difference method for the free boundary value 
problem under logarithmic front-fixing transformation.
In \cite{Fazio:2015:PEE,Fazio:2020:APO}, Fazio proposed one \textit{a posteriori} error estimator for the numerical solution for the American put option 
obtained by an explicit finite difference scheme.

In this paper, we present an implicit finite difference scheme combined with the use of the front-fixing method 
in order to solve the American put option problem without dividend payments. 
Preliminary results have been presented to the International Conference ICNAAM2015 \cite{Fazio:2016:FFF}.
We use a front-fixing variable transformation to reformulate the variable domain problem into a non-linear problem
on a fixed rectangular domain. For the non-linear problem, we introduce a suitable choice of truncated boundary 
that allows to impose the asymptotic boundary condition. Then, the original problem is transformed into a new 
non-linear partial differential equation where the free boundary appears as a new variable and computed as part of 
the solution. We develop an implicit finite difference method, we investigate the consistency and the stability.
Finally, we choose to validate the obtained numerical results via a mesh refinement and the Richardson's extrapolation
and we report the comparison with numerical methods available in the literature.

\section{The American put options model}
Let us suppose that at time $t$ the price of a given underlying asset is $S$.
We consider here the following mathematical model for the value $ P = P(S, \tau)$ of 
an American put option to sell the asset
\begin{align}
& {\ds \frac{\partial P}{\partial \tau}} = \frac{1}{2} \sigma^2 S^2  
{\ds \frac{\partial^2 P}{\partial S^2}} + r S {\ds \frac{\partial P}{\partial S}} -r P \quad \mbox{on} \quad 0 < \tau \le T\ , \ \ S^*(\tau) < S < \infty \ , \nonumber \\[1ex]
&  P(S, 0) = \max(E-S, 0) \ , \qquad S^*(0) = E \ ,  \nonumber \\[1ex]
& \lim_{S\rightarrow \infty} P(S, \tau) = 0 \ , \label{APO:model} \\
& P(S^*(\tau), \tau) = E - S^*(\tau) \ , \qquad {\ds \frac{\partial P}{\partial S}}(S^*(\tau), \tau) = 
- 1 \ , \nonumber \\
& P(S, \tau) = E-S \ , \quad 0 \le S < S^*(\tau) \ , \nonumber
\end{align}
where $\tau = T-t$ denotes the time to maturity $T$, $S^*(\tau)$ is a free boundary, that is the unknown early exercise boundary, $\sigma$, 
$r$ and $E$ are given constant parameters representing the volatility of the underlying asset, the interest rate and the exercise price of the option, respectively.
Equation (\ref{APO:model}) is known as the Black-Scholes-Merton equation and it has been developed by the three economists Fischer Black, Myron Scholes and Robert Merton in $1973$, \cite{Black:POC:1973} and \cite{Merton:1973:TRO}.

\section{The front-fixing method}
The front-fixing method is widely employed for solving free boundary problems.
The basic idea of the front fixing method is to use a variable change in order to remove 
the free boundary and, then, to transform the original equation into a new non-linear partial 
differential equation on a bounded domain, where the free boundary appears 
as a new unknown of the problem. The main advantage of the front-fixing method is that the free 
boundary is computed directly. 

Then, we first transform the Black-Scholes equation into a new parabolic equation over a bounded domain, 
we introduce a truncated boundary and, then, we use finite difference schemes for the new approximate problem
over a bounded domain. 

According to Wu and Kwok \cite{Wu:1997:FFM}, we consider the dimensionless new variables 
\begin{equation}\label{eq:front-fixing}
x = \ln \frac{S}{S^*(\tau)} \ , \quad S_f(\tau) = \frac{S^*(\tau)}{E}\ , \quad p(x, \tau) = \frac{P(S, \tau)}{E} \ ,
\end{equation}
it follows that $S_f(\tau)$ is mapped on the fixed line $x=0$, $0 \le p(x, \tau) \le 1$ and $0 \le S_f(\tau) \le 1$.  
Then, the put option problem (\ref{APO:model}) can be rewritten as follows
\begin{align}
\label{APO:model:fix}
& {\ds \frac{\partial p}{\partial \tau}} = \frac{1}{2} \sigma^2 
{\ds \frac{\partial^2 p}{\partial x^2}} + \left(r - \frac{\sigma^2}{2}\right) {\ds \frac{\partial p}{\partial x}} + \frac{1}{S_f(\tau)} \frac{dS_f}{d\tau}(\tau) {\ds \frac{\partial p}{\partial x}} - r \ p  \ , \\[1ex]
\label{APO:model:fix1}
&  p(x, 0) = 0 \quad \mbox{for} \quad 0 \le x \ , \qquad S_f(0) = 1 \ ,  \\[1ex]
\label{APO:model:fix2}
& \lim_{x \rightarrow \infty } p(x, \tau) = 0 \ , \\
\label{APO:model:fix3}
& p(0, \tau) = 1 - S_f(\tau) \ , \qquad {\ds \frac{\partial p}{\partial x}}(0, \tau) = 
- S_f(\tau)  \ , 
\end{align}
on the domain defined by $0 < \tau \le T$ and $0 < x < \infty$.

In order to solve numerically the obtained problem (\ref{APO:model:fix}-\ref{APO:model:fix3}), we have to consider 
a finite computational domain. Then, we introduce a truncated boundary value $x_\infty$, which is a suitable large value where it would be convenient to impose the asymptotic boundary condition.
In other words, we replace the asymptotic boundary condition (\ref{APO:model:fix2}) with the side condition
\begin{equation}\label{eq:FBC} 
p(x_\infty, \tau) = 0 \ .
\end{equation}
For the choice of  $x_\infty$ and the accuracy of the related numerical solution, we can refer to the study by Kangro and Nicolaides \cite{Kangro:2000:FFB}.

By setting an integer $J$ and a positive value $\mu$, we define the step-sizes
\begin{equation*}\label{eq:step-sizes} 
\Delta x=\frac{x_\infty}{J} \ , \quad \Delta \tau = {\mu}{\Delta x}^{2} \ ,
\end{equation*}
the integer $N$ 
\begin{equation*}\label{eq:N} 
N = \myceil{\frac{T}{\Delta \tau}} \ ,
\end{equation*}
where $ \lceil\cdot\rceil : \RR^+ \rightarrow \NN $ is the \textit{ceiling} function which maps a real number to the least integer that is greater than or equal to that number.
Therefore, $\mu$ is the grid-ratio 
\begin{equation*}\label{eq:c}
\mu = \frac{\Delta \tau}{(\Delta x)^2} \ .
\end{equation*}
Therefore, within the finite domain, we can introduce a mesh of grid-points 
\begin{equation*}\label{eq:grid} 
x_j = j \Delta x \ , \quad \tau^n = n \Delta \tau \ ,
\end{equation*}
for $j =0, 1, \dots , J$ and $n = 0, 1, \dots, N$. 
We would like to define a numerical scheme that allows us to compute the grid values
\begin{equation*}\label{eq:p}
p_j^n \approx p(x_j,\tau^n) \ , 
\end{equation*}
for $j =0, 1, \dots , J$, $n = 0, 1, \dots, N$ and the free boundary values
\begin{equation*}\label{eq:R} 
S_f^n \approx S_f(\tau^n) \ ,
\end{equation*}
for $n = 0, 1, \dots, N$.

\section{A new implicit finite difference scheme}
Here we present an implicit finite difference scheme for (\ref{APO:model:fix}-\ref{APO:model:fix3}).

\begin{align}\label{eq:MN2}
\frac{p_j^{n+1}-p_j^n}{\Delta \tau} &= \frac{1}{2} \sigma^2 \frac{p_{j-1}^{n+1}-2 p_j^{n+1}+p_{j+1}^{n+1}}{(\Delta x)^2} + \\
&+ \left(r - \frac{\sigma^2}{2}\right)\frac{p_{j+1}^{n+1}-p_{j-1}^{n+1}}{2 \Delta x} + \frac{1}{ S_f^{n+1}} \frac{S_f^{n+1}-S_f^{n}}{\Delta \tau}\frac{p_{j+1}^{n+1}-p_{j-1}^{n+1}}{2 \Delta x}-r p_j^{n+1} \nonumber
\end{align}
for $j=1,2,\dots,J-1$ and $n=0,1,\dots,N-1$.\\
The initial conditions (\ref{APO:model:fix1}) are 
\begin{equation*}
p_j^0 = 0 \ , \quad S_f^0 = 1 \ ,
\end{equation*}
for $j= 0,1, \dots , J$.
From the boundary conditions (\ref{eq:FBC}), we get
\begin{equation*}
p_{J}^{n+1} = 0  \ ,
\end{equation*}
for $n= 0,1, \dots , N$.

Moreover, by considering the governing differential equation (\ref{APO:model:fix}) at $x_0 = 0$, $\tau > 0$ and taking into account the side conditions (\ref{APO:model:fix3}), one gets a new boundary condition
\begin{equation}\label{APO:model:fix5}
\frac{\sigma^2}{2}\frac{\partial^2 p}{\partial x^2}(0^+, \tau) + \frac{\sigma^2}{2} S_f(\tau) -r = 0 \ ,
\end{equation}
(see Wu and Kwok \cite{Wu:1997:FFM}, Zhang and Zhu \cite{Zhang:2009:HFD} or Kwok \cite[p. 341]{Kwok:2008:MMF}), and its central finite difference discretization
\begin{equation}\label{eq:FD:fix5}
\frac{\sigma^2}{2}\frac{p_{-1}^{n+1}-2 p_{0}^{n}+p_{1}^{n+1}}{(\Delta x)^2} + \frac{\sigma^2}{2} S_f^{n+1} -r = 0 \ .
\end{equation}
From the two boundary conditions (\ref{APO:model:fix3}), using a central finite difference formula at time $n+1$ and considering (\ref{eq:FD:fix5}), we obtain
\begin {equation} \label{eq:bc:imp}
p_0^{n+1}=1-S_f^{n+1}, \quad p_1^{n+1}= 1+ \frac{r (\Delta x)^2}{\sigma ^2}-\left(1+\Delta x+\frac{(\Delta x)^2}{2}\right)S_f^{n+1}
\end{equation}
where $x_{-1} = -\Delta x$ is a fictitious point out of the computational domain.
Considering $\mu={\Delta \tau}/{(\Delta x)^2}$ and rearranging the (\ref{eq:MN2}), our implicit numerical scheme can be written as
follows
\begin{equation}\label{eq:FDI}
\bar{a}^{n+1} p_{j-1}^{n+1}+\bar{b} p_{j}^{n+1}+\bar{c}^{n+1} p_{j+1}^{n+1}=p_{j}^{n}
\end{equation}
for $j=1,2, \dots,J-1$ and $n=0,1,\dots,N-1$, where
\begin{align}\label{eq:FD:abcMN2}
\bar{a}^{n+1} &= \frac{\mu}{2} \left[- \sigma^2 + \left(r - \frac{\sigma^2}{2}\right) \Delta x \right] + 
  \frac{1}{S_f^{n+1}}\frac{S_f^{n+1}-S_f^{n}}{2 \Delta x}\nonumber \\
\bar{b} &= 1 + \mu \sigma^2 + r \Delta \tau \\
\bar{c}^{n+1} &= \frac{\mu}{2} \left[ -\sigma^2 - \left(r - \frac{\sigma^2}{2}\right) \Delta x \right] -
  \frac{1}{S_f^{n+1}}\frac{S_f^{n+1}-S_f^{n}}{2 \Delta x}\nonumber \ .
\end{align}
Considering (\ref{eq:bc:imp}) and putting $j=1$ in (\ref{eq:FDI}), we get
\begin{equation}\label{eq:j1}
\bar{c}^{n+1}p_2^{n+1}= p_1^{n}-\bar{a}^{n+1}(1-S_f^{n+1})-\bar{b}\left[1+r \frac{(\Delta x)^2}{\sigma^2}-
\left(1+\Delta x+\frac{1}{2} (\Delta x)^2 \right)S_f^{n+1}\right]
\end{equation}
For $j=2$ in (\ref{eq:FDI}) we obtain
\begin{equation}\label{eq:j2}
\bar{b} p_2^{n+1}+\bar{c}^{n+1}p_3^{n+1}= p_2^{n}-\bar{a}^{n+1}\left[1+r \frac{(\Delta x)^2}{\sigma^2}-
\left(1+\Delta x+\frac{1}{2}(\Delta x)^2 \right)S_f^{n+1}\right] \ .
\end{equation}
Putting $j=J$ in (\ref{eq:FDI})
\begin{equation}\label{eq:jJ}
\bar{a}^{n+1}p_{J-1}^{n+1}+\bar{b} p_J^{n+1}= p_J^{n} \ .
\end{equation}
For $j=3, 4, \cdots, J-1$ we have the equations
\begin{equation}\label{eq:j3}
\bar{a}^{n+1}p_{j-1}^{n+1}+\bar{b} p_j^{n+1}+\bar{c}^{n+1}p_{j+1}^{n+1}= p_j^{n} \ .
\end{equation}
Then, at each time step we obtain a system $J$ equations,
given by (\ref{eq:j1})-(\ref{eq:j3}), in $J$ unknowns, $p_2^{n+1},p_3^{n+1},\cdots,p_J^{n+1}$ and $S_f^{n+1}$.
The system (\ref{eq:j1})-(\ref{eq:j3}) can be written in the following compact form 
\begin{equation}\label{eq:sys0}
A(S_f^{n+1})p^{n+1}=f(S_f^{n+1}) \ ,
\end{equation}
where $p^{n+1}=(p_2^{n+1},p_3^{n+1}, \cdots , p_J^{n+1})$, the coefficients matrix $A=A(S_f^{n+1}) \in \mathbb{R}^{J,J-1}$ 
has following form
\begin{eqnarray*}
A(S_f^{n+1})= \left( \begin{array}{cccccc} 
\bar{c}^{n+1} & & & & & 0 \\
\bar{b} & \bar{c}^{n+1} \\
\bar{a}^{n+1}& \bar{b} & \bar{c}^{n+1} \\
& \ddots & \ddots & \ddots &  \\
& &\ddots & \ddots & \ddots &  \\
& & & \bar{a}^{n+1}& \bar{b} & \bar{c}^{n+1}  \\
0 & & & & \bar{a}^{n+1} & \bar{b} \end{array} \right) \ , 
\end{eqnarray*}
and the mapping $f=f(S_f^{n+1}):\mathbb{R}\rightarrow \mathbb{R}^J $ is
\begin{eqnarray*}\label{mat:f}
f(S_f^{n+1})= \left( \begin{array}{l} 
p_1^{n}-\bar{a}^{n+1}(1-S_f^{n+1})-\bar{b}\left[1+r \frac{(\Delta x)^2}{\sigma^2}-
\left(1+\Delta x+\frac{1}{2}(\Delta x)^2 \right)S_f^{n+1}\right] \\
p_2^{n}-\bar{a}^{n+1}\left[1+r \frac{(\Delta x)^2}{\sigma^2}-
\left(1+\Delta x+\frac{1}{2}(\Delta x)^2 \right)S_f^{n+1}\right] \\
p_3^{n} \\
\vdots\\
p_J^{n}  \end{array} \right) .
\end{eqnarray*}
The system (\ref{eq:sys0}) can now be written as a non-linear problem in the form
\begin{equation}\label{eq:sys}
F(p^{n+1},S_f^{n+1})=A(S_f^{n+1})p^{n+1}-f(S_f^{n+1})=0 \ .
\end{equation}
The implicit discretization leads to a system of non-linear equations (\ref{eq:sys})
for the price and the location of the free boundary at each time step.
 
\section{Consistency and stability}
In this section, we discuss the consistency and stability of the implicit finite difference scheme (\ref{eq:MN2}).

\noindent
\subsection{Consistency}
We write the PDE (\ref{APO:model:fix}) and the numerical scheme (\ref{eq:MN2}) as follows
\begin{eqnarray*}
&& L(p, S_f)  =  \frac{\partial p}{\partial \tau} - \frac{1}{2} \sigma ^2 \frac{\partial ^2 p}{\partial x^2}-
\left(r-\frac{\sigma ^2}{2} \right) \frac{\partial p}{\partial x} - \frac{1}{S_f}\frac{\partial S_f}{\partial \tau} \frac{\partial p}{\partial x} +rp=0 \ , 
\end{eqnarray*}
\begin{eqnarray*}
&& F(p^{n+1}_j,S^{n+1}_f)  = \frac{p_j^{n+1}-p_j^n}{\Delta \tau}- \frac{1}{2} \sigma^2 \frac{p_{j-1}^{n+1}-2 p_j^{n+1}+p_{j+1}^{n+1}}{(\Delta x)^2} + \\
&& - \left(r - \frac{\sigma^2}{2}\right)\frac{p_{j+1}^{n+1}-p_{j-1}^{n+1}}{2 \Delta x} 
 - \frac{1}{S_f^{n+1}} \frac{S_f^{n+1}-S_f^{n}}{\Delta \tau}\frac{p_{j+1}^{n+1}-p_{j-1}^{n+1}}{2 \Delta x} + r p_j^{n+1} = 0 .\nonumber
\end{eqnarray*}
In order to study the consistency, let us take an arbitrary point $ (x,\tau) $ in the domain $[0,x_{\infty}] \times [0, T]$, the
mesh point $(x_j, \tau^{n+1})$, a numerical finite difference method is consistent with the differential equation if the local truncation error $T_j^{n+1}$, defined by
$$T_j^{n+1}(\bar{p}, \bar{S}_f)=  F(\bar{p}_j^{n+1}, \bar{S}_f^{n+1})  \ ,$$
 satisfies 
$$T_j^{n+1}(\bar{p},\bar{S}_f) \longrightarrow 0  \qquad as \quad \Delta x \rightarrow 0 \ , \Delta \tau \rightarrow 0  \ , $$
where we denote with $\bar{p}_j^{n+1}=p(x_j,\tau^{n+1})$ the exact solution value of the PDE and with  $\bar{S}_f^{n+1}=S_f(\tau^{n+1})$
the exact solution of the free boundary.
Using Taylor's expansion, assuming the existence of the continuous partial derivatives up to order two in time and up to order four in space, we obtain
\begin{align}\label{ELT}
& T_j^{n+1}(\bar{p},\bar{S}_f) =
 - \frac{\Delta \tau}{2} \ \frac{\partial^2 p}{\partial \tau^2}(x_j,\tau^{n+1}) + O (\Delta \tau)^2
 - \frac{\sigma^2}{2}\frac {(\Delta x)^2}{12} \ \frac{\partial^4 p}{\partial x^4}(x_j,\tau^{n+1}) + O (\Delta x)^4\nonumber \\ 
& - \left(r-\frac{\sigma^2}{2}\right)  \frac{(\Delta x)^2}{3!} \ \frac{\partial^3 p}{\partial x^3}(x_j,\tau^{n+1}) +  O (\Delta x)^4 \nonumber \\ 
& - \frac{1}{S_f^{n+1}} \frac{d S_f}{d\tau}(\tau^{n+1})\left( \frac{(\Delta x)^2}{3!} \ \frac{\partial^3 p}{\partial x^3}(x_j,\tau^{n+1})  + O (\Delta x)^4\right)  \\  
& + \frac{1}{S_f^{n+1}}  \left(\frac{\Delta \tau}{2} \frac{d^2 S_f}{d \tau^2}(\tau^{n+1})\ + O (\Delta \tau)^2 \right)
\left(\frac{\partial p}{\partial x}(x_j,\tau^{n+1}) + \frac{(\Delta x)^2 }{3!}\ \frac{\partial^3 p}{\partial x^3} (x_j,\tau^{n+1}) 
+  O (\Delta x)^4\right) \nonumber \ .
\end{align}
Then, the local truncation error is
\begin{equation}
T_j^{n+1} (\bar{p},\bar{S}_f) = O (\Delta \tau) + O (\Delta x)^2  \ .
\end{equation}
In our case, we have also to consider the boundary conditions (\ref{APO:model:fix3}) and (\ref{APO:model:fix5})
\begin{align*}
L_1(p,S_f)= & p(0,\tau) - 1 +S_f(\tau)=0 \ , \nonumber \\
L_2(p,S_f)= & \frac{\partial p}{\partial x}(0,\tau) +S_f(\tau)=0 \ , \\
L_3(p,S_f)= & \frac{\sigma^2}{2}\frac{\partial^2 p}{\partial x^2}(0,\tau) +\frac{\sigma^2}{2}S_f(\tau) - r=0 \ . \nonumber
\end{align*}
The numerical scheme for the boundary conditions is
\begin{eqnarray*}
&& F_1(p_0^{n+1},S_f^{n+1})=  p_0^{n+1}-1+S_f^{n+1} = 0 \ , \nonumber \\
&& F_2(p_0^{n+1},S_f^{n+1})= \frac{p_1^{n+1}-p_{-1}^{n+1}}{2 \Delta x}+S_f^{n+1} = 0 \\
&& F_3(p_0^{n+1},S_f^{n+1})=  \frac{\sigma^2}{2}\frac{p_{-1}^{n+1}-2p_{0}^{n+1}+p_{1}^{n+1}}{(\Delta x)^2}+ \frac{\sigma^2}{2} S_f^{n+1} -r = 0 \ . \nonumber
\end{eqnarray*}
Using Taylor's expansion, the local truncation error for the boundary conditions is of $O (\Delta x)^2$
\begin{eqnarray*}
&& T_1(\bar{p},\bar{S}_f) =  F_1(\bar{p},\bar{S}_f) = 0 \\
&& T_i(\bar{p},\bar{S}_f) =  F_i(\bar{p},\bar{S}_f) = O(\Delta x)^2 \qquad \text{for} \quad i=2,3 \ .
\end{eqnarray*}

So, assuming the existence of the continuous partial derivatives up to order two in time and up to order four in space
of the solution of the problem (\ref{APO:model:fix}-\ref{APO:model:fix3}), the implicit finite difference method defined by (\ref{eq:MN2}) is consistent with the fixed domain model (\ref{APO:model:fix}) of order $O(\Delta \tau) + O(\Delta x)^2$.

\subsection{Stability} 
Now, we perform the Von Neumann analysis to investigate the stability of the implicit method.
The Von Neumann analysis is only valid for linear problems with constant coefficients, it is not possible to 
apply it to non-linear problems or to problems with variable coefficients.  
In order to apply the method of Von Neumann, it is necessary to linearize the model and freeze the coefficients, 
considering the problem locally \cite{Aslak:IPD:1998}.  Then, the von Neumann analysis can be applied, 
and a stability condition can be derived, this condition will depend on the frozen coefficients involved. 

In our context, the non-linear nature of the model (\ref{APO:model:fix}) is due to the presence of the term
\begin{eqnarray}\label{nl_term}
 \displaystyle{ \frac{1}{S_f}\frac{d S_f}{d \tau}}. 
\end{eqnarray}
with ${S_f}$ unknown of the problem. Then, in the stability analysis, we decide to take it into account and to
derive the stability condition in relation to the value of $S_f$ and of its first derivative. Moreover, note that 
we approximate the term (\ref{nl_term}) with
\begin{eqnarray}\label{S_f}
\frac{1}{S_f^{n+1}} \frac{S_f^{n+1}-S_f^{n}}{ \Delta \tau}  \ .
\end{eqnarray}
that assumes nonpositive values because $S_f^n$ is positive and not increasing.

Using the Fourier analysis, we set
$$ p_j^{n+1}=\lambda p_j^n \qquad p_{j\pm1}^{n+1}=\lambda e^{\pm i k \Delta x}p_j^n$$
where $\lambda=\lambda(k \Delta x)$ is called the amplification factor. 
Substituting these expressions in the numerical method (\ref{eq:FDI}) and dividing by $p_j^n$,  we obtain
$$ \lambda = \frac{1}{ \bar{a}^{n+1} e^{-ik \Delta x} + \bar{b} + \bar{c}^{n+1} e^{ik \Delta x} } \ .$$
By applying the Euler's formulas and by (\ref{eq:FD:abcMN2}), after some manipulations, we derive
\begin{eqnarray*}
\lambda & = & \frac {1}{1+r\Delta \tau+2 \mu \sigma^2 \sin^2 {\displaystyle \frac{k \Delta x}{2}}- i \mu \Delta x \left[\left(r- 
{\displaystyle \frac{\sigma^2}{2}}\right)+ {\displaystyle \frac{1}{ S_f^{n+1}}}
{\displaystyle \frac{S_f^{n+1}-S_f^{n}}{\Delta \tau}} \right]\sin (k \Delta x)} \ .
\end{eqnarray*}
Now if we compute the magnitude of amplification factor $\lambda$, then we have
\begin{eqnarray*}
|\lambda|^2  = \frac{1}{(1 + A)^2 +B^2} \ , 
\end{eqnarray*}
where
\begin{eqnarray*}
&& A =  r\Delta \tau+2 \mu \sigma^2 \sin ^2 {\displaystyle \frac{k \Delta x}{2}} \\
&& B=  \mu \Delta x \left[ \left(r- {\displaystyle \frac{\sigma^2}{2}} \right)
+{\displaystyle \frac{1}{ S_f^{n+1}}} {\displaystyle  \frac{S_f^{n+1}-S_f^{n}}{\Delta \tau }} \right]\sin (k \Delta x) \ . 
\end{eqnarray*}
Because $A>0$, $(1+A)^2> 1$ and $B^2> 0$, thus $|\lambda|<1$.
So, we can conclude that the implicit finite difference method defined by (\ref{eq:MN2}) for the fixed domain problem (\ref{APO:model:fix}-\ref{APO:model:fix3}) is unconditionally stable. 

\begin{figure}[!hbt]
\centering
\psfrag{Phase}[c][][0.8]{ $k \Delta x$}
\psfrag{Amp}[c][][0.8]{$|\lambda|$}
\psfrag{n=1spa}[c][][0.7]{$n=1$}
\psfrag{n=5spa}[c][][0.7]{$n=5$}
\psfrag{n=Nspa}[c][][0.7]{$n=N$}
\includegraphics[scale=.58]{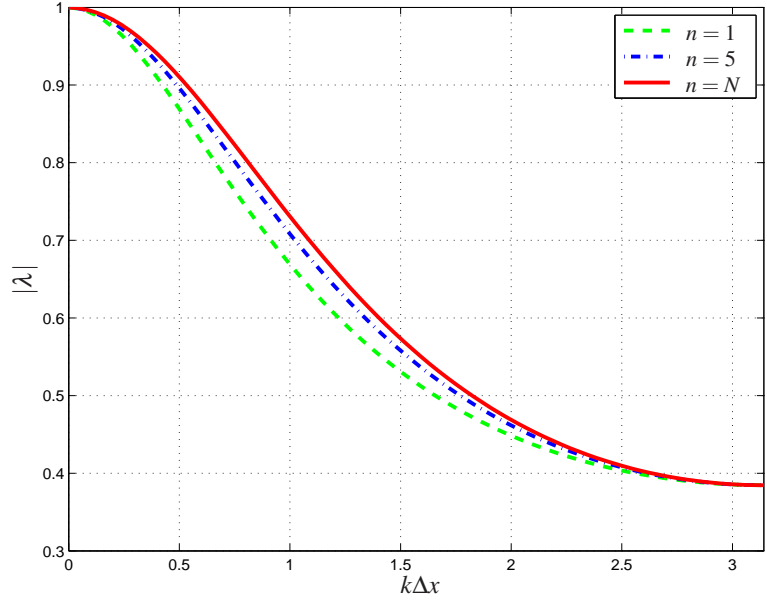}
\caption{\it Implicit method. Amplification factor module $|\lambda|$ for different values of $n$ with $N = 320$, $\mu =20$ and $x_\infty = 1$.}
\label{fig:lambda_n}
\end{figure}

In figure (\ref{fig:lambda_n}), we show the amplification factor module for implicit method at different time steps $\tau^n$. 
In fact, in the stability analysis, we have decided to evaluate the amplification factor $\lambda$ depending
on the value of (\ref{S_f}).
We observe that the variation of (\ref{S_f}) does not influence the stability of the implicit method.

\begin{figure}[!hb]
\centering
\psfrag{Phase}[c][][0.8]{ $k \Delta x$}
\psfrag{Amp}[c][][0.8]{$|\lambda|$}
\psfrag{mu=12spazN=534}[c][][.7]{$\ \mu=12 \quad N=534$}
\psfrag{mu=20spazN=320}[c][][.7]{$\ \mu=20 \quad N=320$}
\psfrag{mu=26spazN=247}[c][][.7]{$\ \mu=26 \quad N=247$}
\psfrag{mu=100spaN=64}[cc][][.7]{$\mu=100 \ \  N=64$}
\includegraphics[scale=.55]{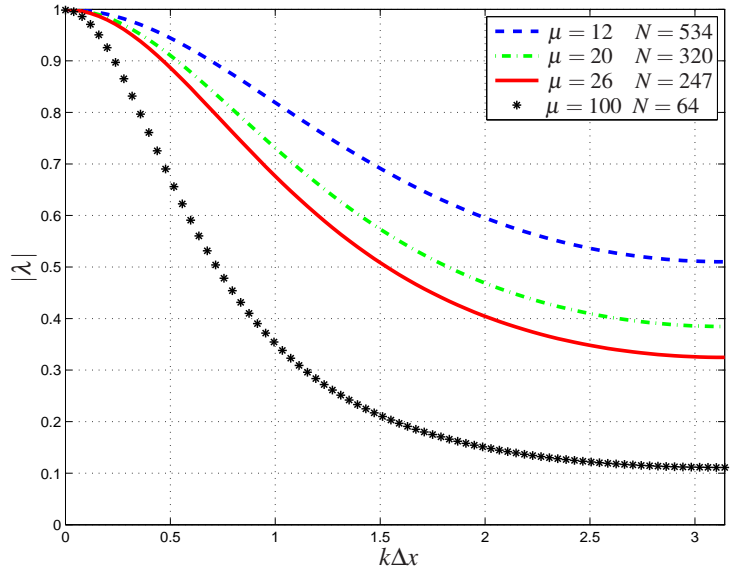}
\caption{\it Implicit method. Amplification factor module $|\lambda|$ for different values of $\mu$ at time $\tau=T$. 
The implicit method results to be stable for any value of $\mu$.}
\label{fig:lambda_imp}
\end{figure}

In figure (\ref{fig:lambda_imp}) we show the amplification factor module for different values of $\mu$ at final time $\tau=T$. 
We observe that the unconditionally stable implicit finite difference scheme produces results
that are stable for any value of grid-ratio $\mu$.

\begin{figure}[!ht]
\centering
\psfrag{Phase}[c][][0.8]{ $k \Delta x$}
\psfrag{mu=12spazN=534}[c][][0.7]{$\mu=12 \quad N=534$}
\psfrag{mu=20spazN=320}[c][][0.7]{$\mu=20 \quad N=320$}
\psfrag{mu=26spazN=247}[c][][0.7]{$\mu=26 \quad N=247$}
\psfrag{Ampl}[c][][0.8]{$|\lambda|$}
\includegraphics[scale=.6]{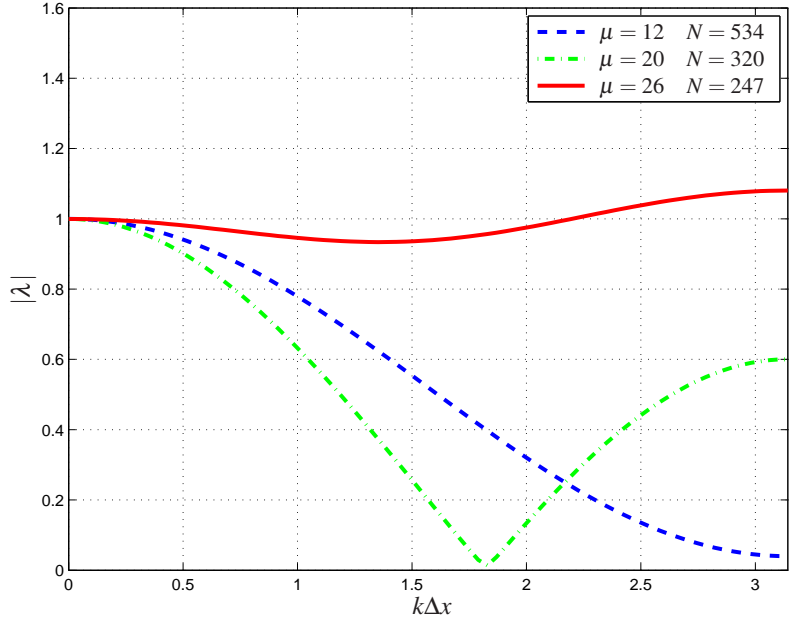}
\caption{\it Explicit method. Amplification factor module $|\lambda|$ for different values of $\mu$ at time $\tau=T$. The explicit method
becomes unstable for $\mu \geq 26$.}
\label{fig:lambda_exp}
\end{figure}

We can compare the result reported in figure (\ref{fig:lambda_imp}) with the ones reported in figure (\ref{fig:lambda_exp}), where
we show the amplification factor module $\lambda$, for different values of $\mu$, of the explicit method \cite{Company:2014:SAP}
\begin{eqnarray*}
\lambda & = & {1-r\Delta \tau - 2 \mu \sigma^2 \sin^2 {\displaystyle \frac{k \Delta x}{2}} + i \mu \Delta x \left[\left(r- 
{\displaystyle \frac{\sigma^2}{2}}\right)+ {\displaystyle \frac{1}{ S_f^{n}}}
{\displaystyle \frac{S_f^{n+1}-S_f^{n}}{\Delta \tau}} \right]\sin (k \Delta x)} \ .
\end{eqnarray*}
In order to evaluate $\lambda$, we use the same assumptions given for the implicit method.
We observe that the explicit finite difference scheme produces results
that are stable only for values of grid-ratio $\mu < 26$.

The main limitation of the explicit method is the restriction on the choice of $\Delta \tau$ and $\Delta x$ that have to be
chosen in relation at the values of the parameters of the model $\sigma$ and $r$, the volatility of the underlying asset and the interest rate, respectively. On the other hand, the new implicit method is unconditionally stable, no restriction on step sizes is required.

\newpage
\section{\textit{A posteriori} error estimator}
In this section, we describe \textit{a posteriori} error estimator for the American put options problem based on Richardson's extrapolation. In this way, we can find a more suitable computational grid requiring that the numerical solution obtained by the implicit method verifies a prefixed error tolerance. 

For a scalar $U$ of interest, either a value of the solution $p_j^n$ or a free boundary value $S_f^n$, the numerical error $e$ can be defined by
\begin{equation}\label{eq:GE}
e = u - U \ ,
\end{equation}
where $u$ is the exact, usually unknown, value.
When the numerical error is caused prevalently by the discretization error and in the case of smooth enough solutions the global error can be split into a sum of powers of the inverse of $N$ 
\begin{equation}\label{eq:asymE}
u = U_N + C_0 \left(\frac{1}{N}\right)^{q_0}+ C_1 \left(\frac{1}{N}\right)^{q_1}+ C_2 \left(\frac{1}{N}\right)^{q_2}+ \cdots \ ,
\end{equation}
where $C_0$, $C_1$, $C_2$, $\dots$ are coefficients that depend on $u$ and its derivatives, but are independent of $N$, and $q_0$, $q_1$, $q_2$, $\dots$ are the true orders of the discretization error, see Schneider and Marchi \cite{Schneider:GRR:2005} and the references quoted therein.
Each $q_k$ is usually a positive integer with $q_0 < q_1 < q_2 < \cdots$ and all together constitute an arithmetic progression of ratio $q_1-q_0$.
The value of $q_0$ is called the asymptotic order or the order of accuracy of the method or of the numerical value $U_N$. 
By replacing into equation (\ref{eq:asymE}) $N = N_g$ and $N = N_{g+1}$ and subtracting, to the second obtained equation the first times $(1/s)^{q_0}$, $s = N_{g+1}/N_{g}$, we get the first extrapolation formula 
\begin{equation}\label{eq:Rextra1}
u \approx  U_{g+1} + \frac{U_{g+1}-U_{g}}{s^{q_0}-1} \ ,
\end{equation}
that has a leading order of accuracy equal to $q_1$.
This type of extrapolation is due to Richardson \cite{Richardson:1910:DAL,Richardson:1927:DAL}.
Taking into account equation (\ref{eq:Rextra1}), we can conclude that the error estimator by a first Richardson's extrapolation is given by
\begin{equation}\label{eq:est1}
e_{r} = \frac{U_{g+1}-U_{g}}{s^{q_0}-1} \ ,
\end{equation}
where $q_0$ is the order of the numerical method used to compute the numerical solutions.
Hence, (\ref{eq:est1}) gives the real value of the numerical solution error without knowledge of the exact solution.
In comparison with (\ref{eq:est1}), a safer error estimator can be defined by
\begin{equation}\label{eq:est2}
e_{s} = U_{g+1}-U_{g} \ .
\end{equation}
Of course, $q_0$ can be estimated with the formula
\begin{equation}\label{eq:p0}
q_0 \approx {\ds \frac{\log(|U_g-u|)-\log(|U_{g+1}-u|)}{\log(s)}} \ ,
\end{equation}
where $u$ is again the exact solution (or, if the exact solution is unknown, a reference solution computed with a suitable large value of $N$), and both $u$ and $U_{g+1}$ are evaluated at the same grid-points of $U_g$. 

Within the above framework, in order to improve the numerical accuracy by using only a small number of grid-nodes, we can generalize (\ref{eq:Rextra1}) introducing the following repeated extrapolation formula
\begin{equation}\label{eq:Rextra}
U_{g+1,k+1} = U_{g+1,k} + \frac{U_{g+1,k}-U_{g,k}}{s^{q_k}-1} \ ,
\end{equation}
where $g \in \{0, 1, 2 , \dots , G-1\}$, $k \in \{0, 1, 2, \dots , G-1\}$, $s = N_{g+1}/N_{g}$ is the grid refinement ratio, and $q_k$ is the true order of the discretization error.
The formula (\ref{eq:Rextra}) is asymptotically exact in the limit as $N_0$ goes to infinity if we use uniform grids.
We notice that to obtain each value of $U_{g+1,k+1}$ we need two computed solutions $U$ in two adjacent grids, namely $g+1$ and $g$ at the extrapolation level $k$.
For any $g$, the level $k=0$ represents the numerical solution of $U$ without any extrapolation.
We recall that the theoretical orders of accuracy of the numerical values $U_{g,k}$, with $N =N_g$ and $k$ extrapolations, verify the relation
\begin{equation}\label{eq:pk}
q_k = q_0 + k (q_1-q_0) 
\end{equation}
valid for $k \in \{0, 1, 2, \dots , G-1\}$.

\section{Numerical results}
In this section, we show numerical results obtained by using the implicit finite difference scheme.
Comparisons with some numerical methods available in the literature for the American put options problem are carried out.

The parameters considered for solving the American put option problem (\ref{APO:model:fix}-\ref{APO:model:fix3}) are the following
\begin{equation}\label{eq:parameters}
r = 0.1 \ , \quad \sigma = 0.2 \ , \quad  T = 1 \ , \quad E = 1\ .
\end{equation}
In order to choose the truncated boundary value $x_\infty$, we computed the numerical solution $S_{f}^{N}$ for three different values of $x_\infty$. As we can see in table \ref{tab:xinf}, the change of these values does not affect much the numerical solution; for this reason, we decided to set $x_\infty=1$.\\
\begin{table}[!hbt]
\centering
\renewcommand\arraystretch{1.2}
\begin{tabular}{cccc}
\hline
{$x_\infty$} & {$N = 10$} & {$N = 20$} & {$N = 40$} \\
\hline  
$1$  & $0.8710513685210$  & $0.8661003514438$ & $0.86351373597835$  \\ 
$2$  & $0.8710513685385$  & $0.8661003514444$ & $0.86351373597828$  \\ 
$4$  & $0.8710513685388$  & $0.8661003514438$ & $0.86351373597789$   \\ 
\hline
\end{tabular}
\caption{Free boundary value $S_f^N$ at $\tau=T$ for different truncated boundary locations $x_{\infty}$.}
\label{tab:xinf}
\end{table}

\begin{figure}[!]
\centering
\psfrag {x}[c][][0.7]{$x$} \psfrag {p}[c][][0.7]{$p_j^N \approx p_j(T)$}
\includegraphics[width=0.49\textwidth]{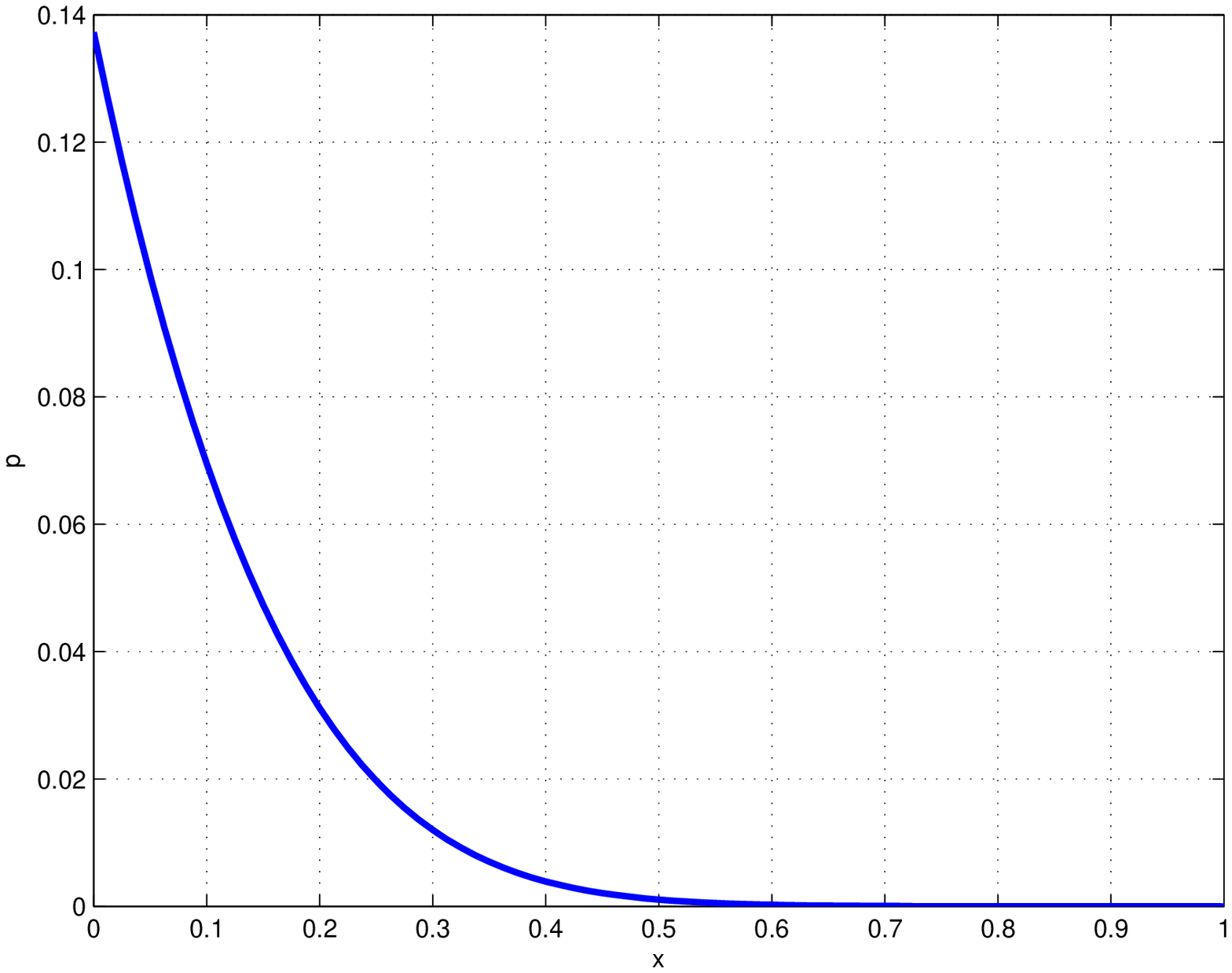} 
\hfill
\psfrag {t}[c][][0.7]{$\tau$} \psfrag {S}[c][][0.7]{${S_f}^n \approx S_f(\tau^n)$}
\includegraphics[width=0.49\textwidth]{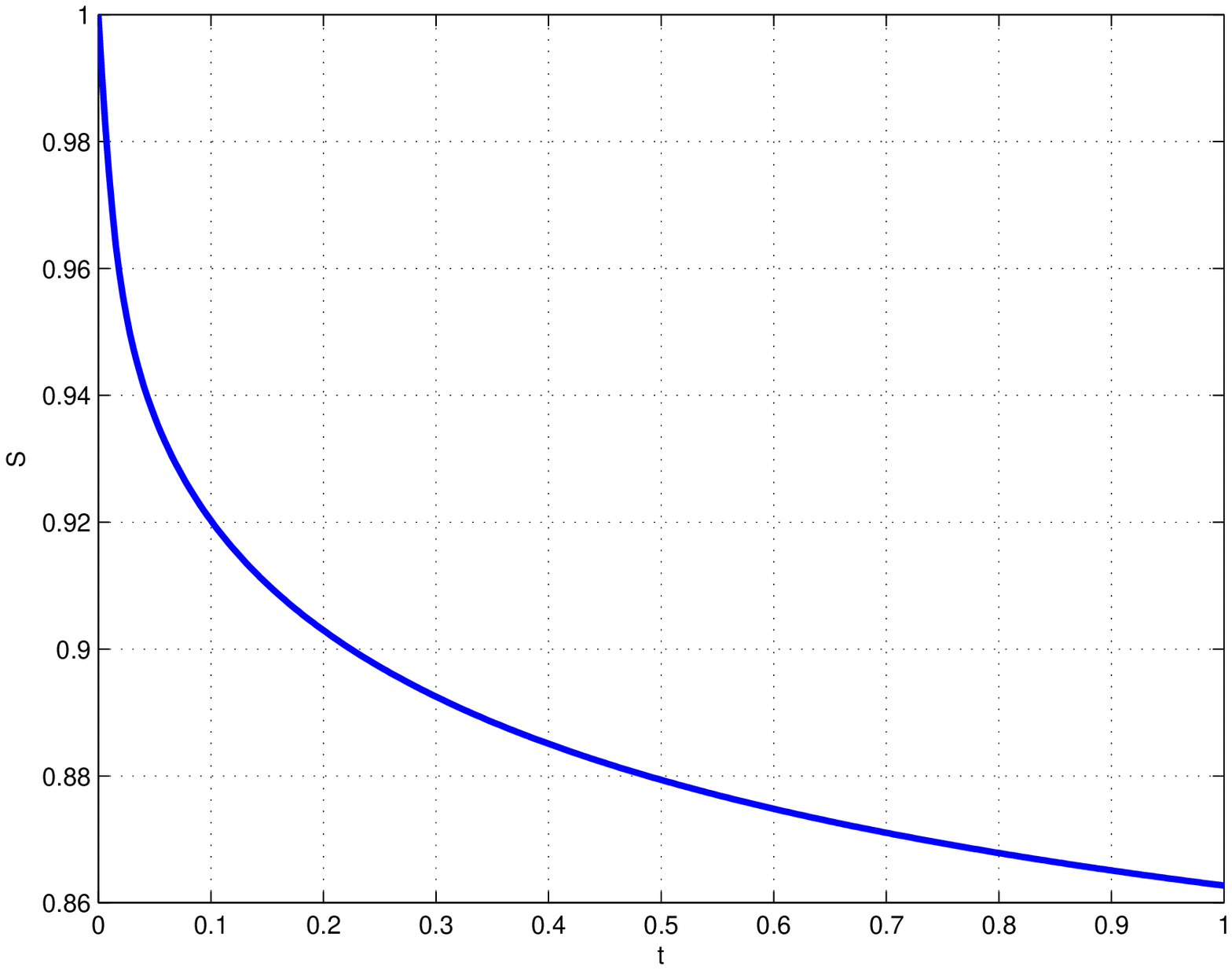}
\caption{\it Numerical results obtained by implicit scheme (\ref{eq:MN2}): on the left  $p_j^N$ versus $x_j$, on the right, $S_f^n$ versus $\tau^n$, obtained with $J=80$, $N=320$ and $\mu=20$.}
\label{fig:impfdAPO}
\end{figure}

\begin{figure}[!]
\centering
\psfrag {x}[c][][0.7]{$x$} \psfrag {p}[c][][0.7]{$p_j^n$}
\psfrag {t}[c][][0.7]{$\tau$}
\includegraphics[width=0.8\textwidth]{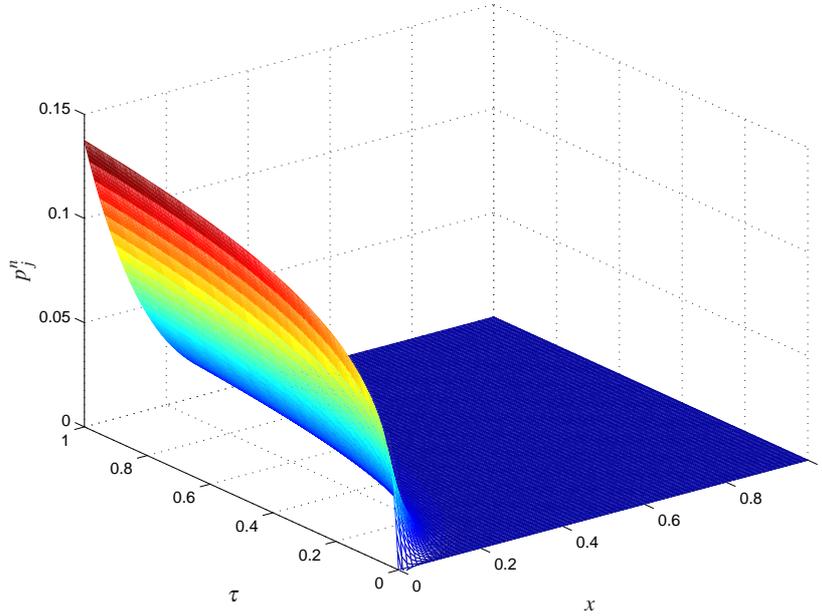} 
\caption{\it Numerical results obtained by implicit scheme (\ref{eq:MN2}), $p_j^n$ with $J=80$, $N=320$ and $\mu=20$.}
\label{fig:impfdAPO_mesh}
\end{figure}

In figure \ref{fig:impfdAPO}, we show the plots of $p_j^N$ versus $x_j$ and of $S_f^n$ versus $\tau^n$; these results are obtained by the implicit finite difference scheme, considering the parameters (\ref{eq:parameters}) with $N = 320$ and $\mu = 20$. 
In figure (\ref{fig:impfdAPO_mesh}), the $3$D plot of numerical solution $p_j^n$ is shown.

In order to solve the non-linear system, we use the MATLAB routine \lq \lq fsolve\rq \rq.
Initially, we used the classical Newton iteration method, but the obtained numerical results showed less accurate 
when compared with ones obtained by Matlab routine \emph{fsolve}.
By our numerical experiments, we conclude that the MATLAB routine, principally regarding the choice of the initial iterate, reveals more robust than the Newton iteration method 
and then more suited for our application.

Looking at the results listed in Table \ref{tab:xinf}, we realize also that, fixing a value of the truncated boundary $x_{\infty}$, the computed values of $S_f^N$ for different values of the grid-steps are in agreement only up to the first decimal place. For this reason, we decide to improve the numerical accuracy by performing a mesh refinement applying repeated Richardson's extrapolation.
The implicit difference scheme is accurate of first-order in time and second-order in space, then, we can perform a mesh refinement
keeping constant the grid-ratio $\mu$, so that we end up with second-order truncation error $T_j^n = O(\Delta \tau^2)$ in time.
Therefore, the global error is of first-order, which is the value of $q_0 =1$.
In our case the sequence of $s^{q_k}$, for $k=0,1, \dots$, with $s = 4$ and $q_k = k+1$ is given by $4$, $16$, $64$, $256$, $1024$, $\dots$. \\
\begin{table}[!hbt]
\renewcommand\arraystretch{1.2}
\centering
\begin{tabular}{r|cccccc}
\hline
{$N$} & {$U_{g,0}$} & {$U_{g,1}$} & {$U_{g,2}$} & {$U_{g,3}$} & {$U_{g,4}$} & {$U_{g,5}$}\\
\hline  
$5$         & $0.884069$ & $$ & $$ & $$ & $$ \\
$20$        & $0.866100$ & $0.860111$ & $$ & $$ & $$ \\
$80$        & $0.863100$ & $0.862100$ & $0.862232$ & $$ & $$ \\
$320$       & $0.862719$ & $0.862592$ & $0.862625$ & $0.862631$ & $$ \\
$1280$      & $0.862717$ & $0.862716$ & $0.862724$ & $0.862726$ & $0.862726$ & $$\\
$5120$      & $0.862738$ & $0.862746$ & $0.862748$ & $0.862748$ & $0.862748$ & $0.862748$\\
\hline 
\end{tabular}
\caption{Implicit method: Richardson's repeated extrapolations for the free boundary value  $S_f^N$ at $\tau=T$.}
\label{tab:Rextraimp}
\end{table}
The numerical results, obtained by repeated Richardson's extrapolations, for the values $S_f^N$ computed with the implicit method are reported in Table \ref{tab:Rextraimp}. The last extrapolated is $U_{5,5}=0.862748$, so we can consider this value as our benchmark $S_f^{N} \approx 0.862748$.
Our result can be compared with the value $S_f^{N} \approx 0.862762$ found by Fazio \cite{Fazio:2015:PEE}, with $S_f^{N} \approx 0.86222$ computed by Nielsen et al. \cite{Nielsen:2002:PFF} and with $S_f^{N} \approx 0.8628$ found by the explicit method of Company et al. \cite{Company:2014:SAP}.\\
In Table \ref{tab:valuetab}, we compare the option price $P(S,T)$ obtained by different methods with the following parameters 
\begin{eqnarray}\label{param}
T = 3 \qquad \sigma = 0.2 \qquad r = 0.08   \qquad E = 100 \ .
\end{eqnarray}
We report the \lq \lq true value\rq \rq \, as the reference offered in \cite{Saib:2011:SFD}, the penalty method (PM) of Nielsen et al. given in \cite{Saib:2011:SFD}, the explicit method (EM) of Company et al. \cite{Company:2014:SAP} with $\Delta x=0.02$, the explicit method with the Richardson's extrapolation (EMR) of Fazio \cite{Fazio:2015:PEE} with $\mu=5$ and our implicit method (IM) without extrapolation, setting $\Delta x =0.02$ and $\mu=5$, and with Richardson's extrapolation (IMR) with $\mu=5$.
In order to find the value of the option price $P(S,T)$ in correspondence to each
different asset price $S$ we use a piecewise cubic spline interpolation.

\begin{table}[!hbt]
\centering
\renewcommand\arraystretch{1.2}
\begin{tabular}{cccccccc}
\hline
\text{Asset} & \text{True} & \text{PM} & \text{EM} & \text{EMR} & \text{IM} & \text{IMR}\\
   {price S} & {value}     &  Nielsen \cite{Nielsen:2002:PFF} &Company \cite{Company:2014:SAP} &  Fazio \cite{Fazio:2015:PEE}&
	& &\\
\hline  
$90$   & $11.6974$ & $11.7207$ & $11.7054$  & $11.7706$  & $11.6926$ & $11.7707$ \\ 
$100$  & $6.9320$  & $6.9573$  & $6.9309$   & $6.9313$   & $6.9243$  & $6.9315$ \\ 
$110$  & $4.1550$  & $4.1760$  & $4.1564$   & $4.1288$   & $4.1467$  & $4.1291$ \\ 
$120$  & $2.5102$  & $2.5259$  & $2.5151$   & $2.5061$   & $2.5028$  & $2.5064$ \\ 
\hline
\end{tabular}
\caption{Comparison of American put option price $P(S,T)$ with parameters (\ref{param}).}
\label{tab:valuetab}
\end{table}
Our goal is to solve the American put option problem with a given tolerance $\e$, where $0 < \e \ll 1$; then 
we use the error estimator defined by equation (\ref{eq:est1}), or alternatively by equation (\ref{eq:est2}).
To this end, we should solve the given problem twice, for two grids defined with given values of $J_g = J$ and $J_{g+1} =2 J$ of space intervals but for the same value of the grid-ratio $\mu$.
The corresponding time intervals $N_g$ and $N_{g+1}$ verify the relation $s = N_{g+1}/N_{g}$.
Hence, we can apply (component-wise) to $p^n$ and to $S_f^n$ the error estimator formula (\ref{eq:est1}), or (\ref{eq:est2}). 
Then, we can verify whether, for $n = 1, 2, \dots, N$,
\begin{equation}\label{eq:test}
\norme{e_r(p^n)}{\infty} \le \e \qquad |e_r(S_f^n)| \le \e \ .
\end{equation}
If the two inequalities (\ref{eq:test}) hold true, for $n = 1, 2, \dots, N$, then we can accept the numerical solution computed on the grid defined by $J_{g+1}$ and $N_{g+1}$, otherwise, we have to increase these two integers and repeat the computation.

\begin{figure}[!t]
\centering
\psfrag {t}[c][][0.7]{$\tau$} \psfrag {er}[c][][0.7]{$e_r(p_j^n)$, $e_r(S_f^n)$}
\psfrag{errp}[c][][0.5]{$e_r(p_j^n)$}
\psfrag{errs}[c][][0.5]{$e_r(S_f^n)$}
\includegraphics[width=0.49\textwidth]{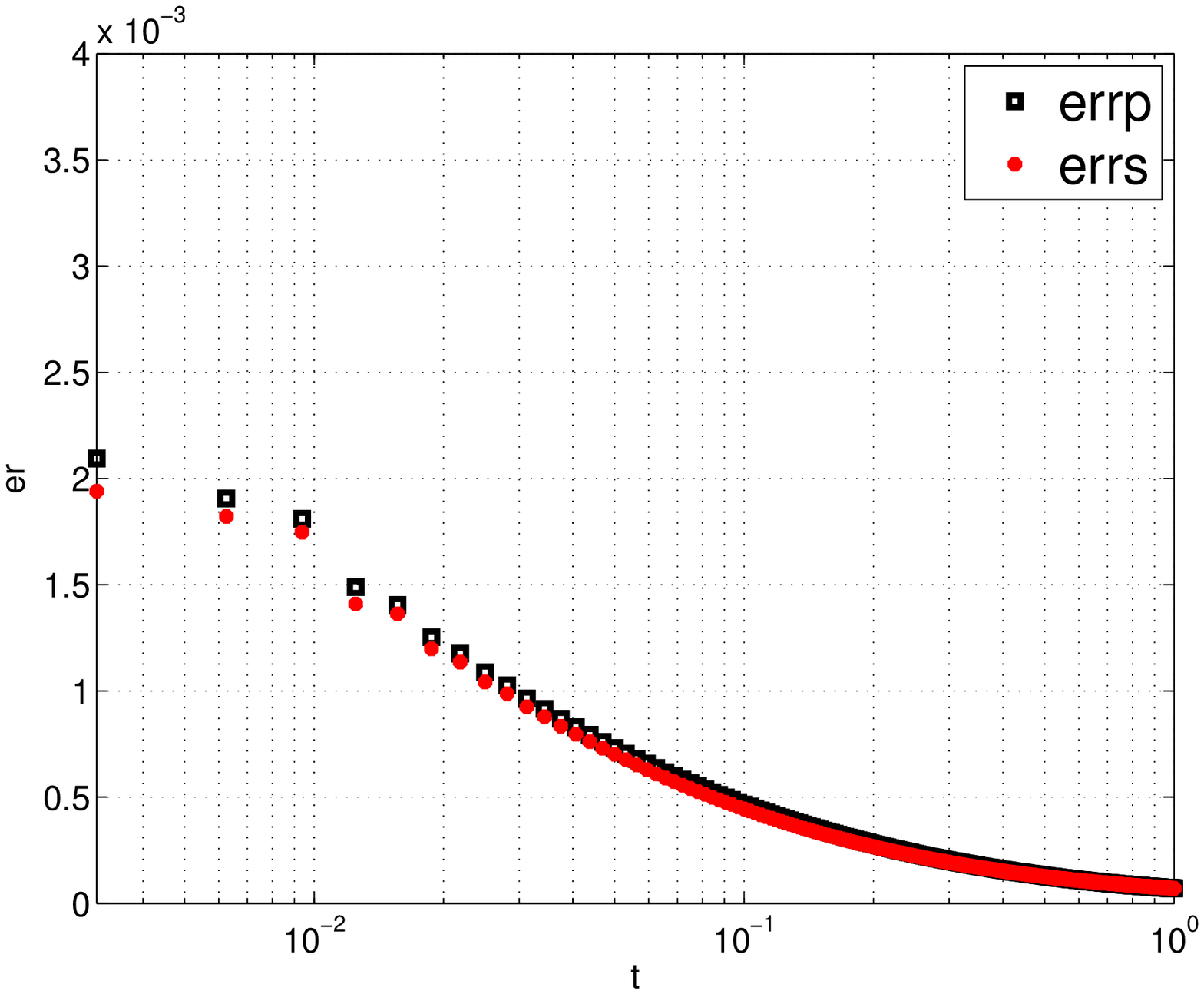}
\hfill
\psfrag {t}[c][][0.7]{$\tau$} \psfrag {er}[c][][0.7]{$e_r(p_j^n)$, $e_r(S_f^n)$}
\psfrag{errp}[c][][0.5]{$e_r(p_j^n)$}
\psfrag{errs}[c][][0.5]{$e_r(S_f^n)$}
\includegraphics[width=0.49\textwidth]{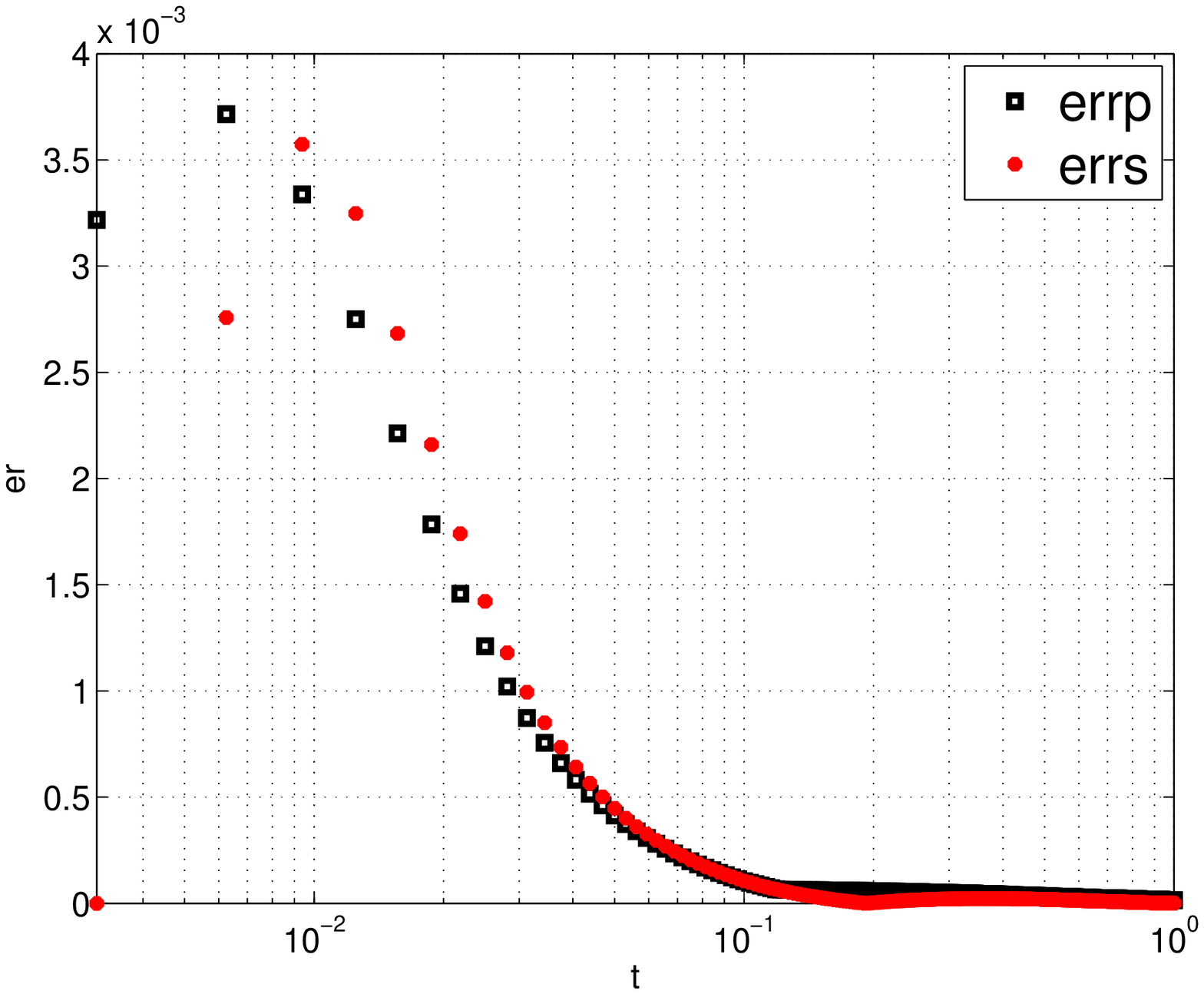}
\caption{\it Numerical estimated errors $e_r(p^n)$ and $e_r(S_f^n)$ versus $\tau^n$ for the explicit (left) and implicit method (right), setting $\epsilon = 0.005$, $\mu=20$, $J_5=160$ and $N_5=1280$.}
\label{fig:APOe}
\end{figure}

Fig. \ref{fig:APOe} shows the error estimator results computed by setting $\e = 0.005$.
We set $\mu=20$ and start with $J_0=5$ and $J_1= 10$ repeating the computation by doubling the number of spatial grid-intervals if the required accuracy is not achieved. 
Our algorithm stops when $J_5 = 160$ that for $\mu = 20$ means $N_5=1280$.
From the numerical results shown in Figure \ref{fig:APOe} we can easily realize that the greatest errors are found within a few time steps
for both numerical methods. 
The error of the implicit method is almost twice of the error of the explicit scheme only in the first time steps, but it becomes smaller in the next time steps, so we can conclude that the numerical results obtained by the implicit method are more accurate than those of 
the explicit one.
We suggest that in order to obtain better accuracy also in the first time steps it can be developed an adaptive version of the both finite difference schemes. 

\section{Concluding Remarks}
In this paper, we consider an American put option model, a free boundary problem defined on an infinite domain.
We overcome the numerical difficulty of solving a free boundary problem using a front fixing approach combined with
finite difference schemes. 
We propose an unconditionally stable implicit finite difference scheme and we prove
the consistency. In order to improve the accuracy solution, we use Richardson's extrapolation.
Comparisons with some methods available in the literature are carried out to validate the obtained numerical results.
Finally, by \textit{a posteriori} error estimator, we find a suitable computational grid requiring that the computed solution
verifies a prefixed tolerance.

\vspace{1.5cm}

\noindent {\bf Acknowledgement.} {The research of this work was 
partially supported by the University of Messina and by the GNCS of INDAM.}

\end{document}